\documentclass[11pt]{amsart}

\usepackage{amsfonts}
\usepackage{amscd}
\usepackage{amsmath, mathrsfs, amssymb}
\usepackage{amsthm}
\usepackage{setspace}
\usepackage{hyperref}
\usepackage{color}
\usepackage{epsfig}
\usepackage{here}
\usepackage{graphicx}
\usepackage[all]{xy}
\usepackage{psfrag}
\usepackage{graphicx,transparent}
\usepackage{enumerate}
\usepackage{caption}
  
\theoremstyle{plain}

\theoremstyle{definition}

\newcommand{\MM}{\mathcal M}

\newcommand{\BM}{\overline{\mathcal M}}

\newcommand{\calN}{\mathcal N}
\newcommand{\calQ}{\mathcal Q}

\newcommand{\calX}{\mathcal X}

\newcommand{\OO}{\mathcal O}

\newcommand{\RR}{\mathbb R}
\newcommand{\BHH}{\overline{\mathcal H}}

\newcommand{\calP}{\mathcal P}

\newcommand{\HH}{\mathcal H}

\newcommand{\Area}{\operatorname{Area}}

\newcommand{\hyp}{\operatorname{hyp}}
\newcommand{\GL}{\operatorname{GL}}
\newcommand{\SO}{\operatorname{SO}}

\newcommand{\odd}{\operatorname{odd}}

\newcommand{\cyl}{\operatorname{cyl}}

\newcommand{\ord}{\operatorname{ord}}

\newcommand{\Hyp}{\operatorname{Hyp}}
\newcommand{\bbC}{\mathbb C}

\newcommand{\bbP}{\mathbb P}

\newcommand{\bbQ}{\mathbb Q}
\newcommand{\bbR}{\mathbb R}
\newcommand{\bbZ}{\mathbb Z}
\newcommand{\bbH}{\mathbb H}

\newcommand{\SL}{\operatorname{SL}}

\newcommand{\Res}{\operatorname{Res}}

\newcommand{\TT}{\mathcal{T}}

\makeatletter
	\@namedef{subjclassname@2010}{%
	\textup{2010} Mathematics Subject Classification}
	\makeatother

\title{Teichm\"uller dynamics in the eyes of an algebraic geometer}

\author{Dawei Chen}
\address{Department of Mathematics, Boston College, Chestnut Hill, MA 02467}
\email{dawei.chen@bc.edu}


\date{}

\thanks{During the preparation of this article the author was partially supported by the NSF CAREER Award DMS-1350396.}

\begin{document}

\begin{abstract}
This is an introduction to the algebraic aspect of Teichm\"uller dynamics, with a focus on its interplay with the geometry of moduli spaces of curves as well as recent advances in the field. 
\end{abstract}

\maketitle

\setcounter{tocdepth}{1}
\tableofcontents

\section{Introduction}
\label{sec:intro}

\begin{figure}[h]
    \centering
    \psfrag{a}{$a$}
    \psfrag{b}{$b$}
    \psfrag{c}{$c$}
    \psfrag{d}{$d$}
    \includegraphics[scale=0.5]{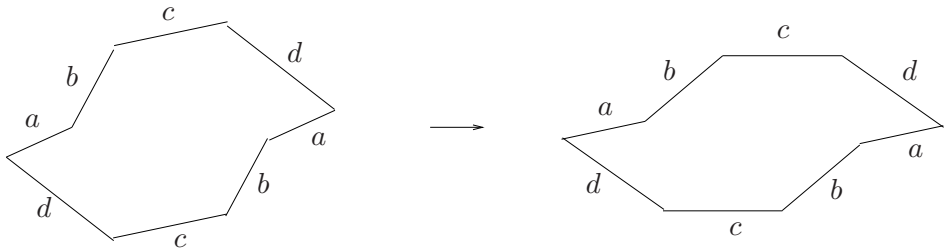}
    \caption{\label{fig:teich} $\GL_2^+(\RR)$-action on a flat surface}
    \end{figure}

An Abelian differential defines a flat structure such that the underlying
Riemann surface can be realized as a plane polygon whose edges are pairwise identified via translation. Varying the shape of the polygon by $\GL_2^+(\RR)$ induces an action on the moduli space of Abelian differentials, called Teichm\"uller dynamics, see Figure~\ref{fig:teich}. 

The corresponding $\GL_2^+(\RR)$-orbit closures in the moduli space of Abelian differentials are now known as affine invariant submanifolds. A number of questions about surface geometry boil down to understanding the structures of affine invariant submanifolds. From the viewpoint of algebraic geometry affine invariant submanifolds are of an independent interest, which can provide special subvarieties in the moduli space of curves. 

We aim to introduce Teichm\"uller dynamics from the viewpoint of algebraic geometry. In Section~\ref{sec:prelim} we review background material, including translation surfaces, strata of Abelian differentials, the $\GL_2^+(\RR)$-action, and affine invariant submanifolds. Section~\ref{sec:teich} focuses on Teichm\"uller curves formed by closed $\GL_2^+(\RR)$-orbits, where we describe their properties, examples, classification, and invariants. In Section~\ref{sec:affine} we study general affine invariant submanifolds and survey recent breakthroughs about their structures, classification, and boundary behavior. Finally in Section~\ref{sec:higher} we discuss similar questions for meromorphic and higher order differentials. 

This article is written in an expository style. We will often highlight motivations and minimize technical details. For further reading, we refer to \cite{Zorich, MoellerSurvey, Wright1, Wright2} for a number of excellent surveys on related topics. 

\subsection*{Acknowledgement}
This article is partially based on the lectures given by the author during the Algebraic Geometry Summer Research Institute Bootcamp, July 2015. The author is grateful to the Bootcamp organizers Izzet Coskun, Tommaso de Fernex, Angela Gibney, and Max Lieblich for their invitation and hospitality. The author thanks Alex Wright for carefully reading the article and many helpful comments.  


\section{Preliminaries}
\label{sec:prelim}

In this section we introduce basic background material that will be used later. 

\subsection{Abelian differentials and translation surfaces}
\label{subsec:abelian}

\begin{figure}[h]
    \centering
    \psfrag{B1}{$B_1$}
    \psfrag{B2}{$B_2$}
    \psfrag{Bk+1}{$B_{k+1}$}
    \psfrag{A1}{$A_1$}
    \psfrag{A2}{$A_2$}
    \psfrag{A3}{$A_3$}
    \psfrag{Ak+1}{$A_{k+1}$}   
    \includegraphics[scale=1.2]{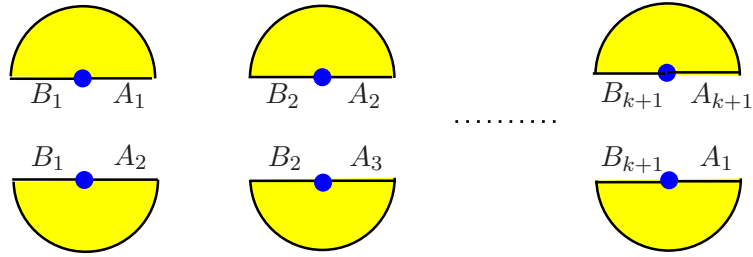}
    \caption{\label{fig:zero-k} A saddle point of cone angle $(k+1)\cdot (2\pi)$}
    \end{figure}

A \emph{translation surface} (also called a \emph{flat surface}) is a closed, topological surface $X$ together with a finite set 
$\Sigma\subset X$ such that: 
\begin{itemize}
\item 
There exists an atlas of charts $X\backslash \Sigma \to \bbC$, where the transition functions are translation. 
\item 
For each $p\in \Sigma$, under the Euclidean metric of $\bbC$ 
 the total angle at $p$ is $(k+1)\cdot (2\pi) $ for some $k\in \bbZ^+$. 
\end{itemize}
We say that $p$ is a \emph{saddle point of cone angle} $(k+1)\cdot (2\pi)$. Locally one can glue $2k+2$ half-disks consecutively to form a cone of angle $2\pi \cdot (k+1)$, see Figure~\ref{fig:zero-k}. 

Equivalently, a translation surface is a closed Riemann surface $X$ with an \emph{Abelian differential} $\omega$, not identically zero: 
\begin{itemize}
\item 
The set of zeros of $\omega$ corresponds to $\Sigma$. 
\item 
If $p$ is a zero of $\omega$ of order $k$, then the cone angle at $p$ is $(k+1)\cdot (2\pi)$. 
\end{itemize}

\begin{figure}[h]
    \centering
    \psfrag{a}{$v_1$}
    \psfrag{b}{$v_2$}
    \psfrag{c}{$v_3$}
    \psfrag{d}{$v_4$}
    \includegraphics[scale=0.5]{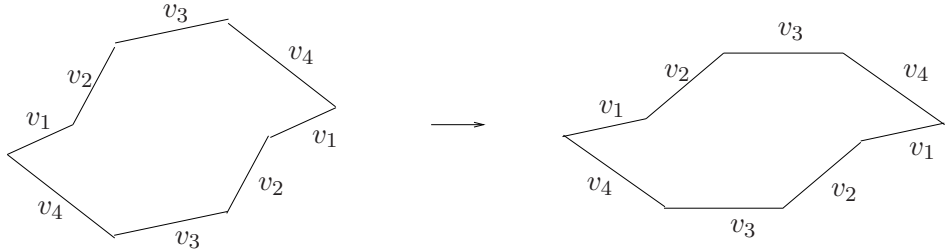}
    \caption{\label{fig:flat} An octagon $X$ with four pairs of parallel edges}
    \end{figure}

\begin{figure}[h]
    \centering
    \psfrag{p}{$p$}
    \psfrag{z}{$z$}
    \psfrag{w}{$z'$}
    \includegraphics[scale=0.5]{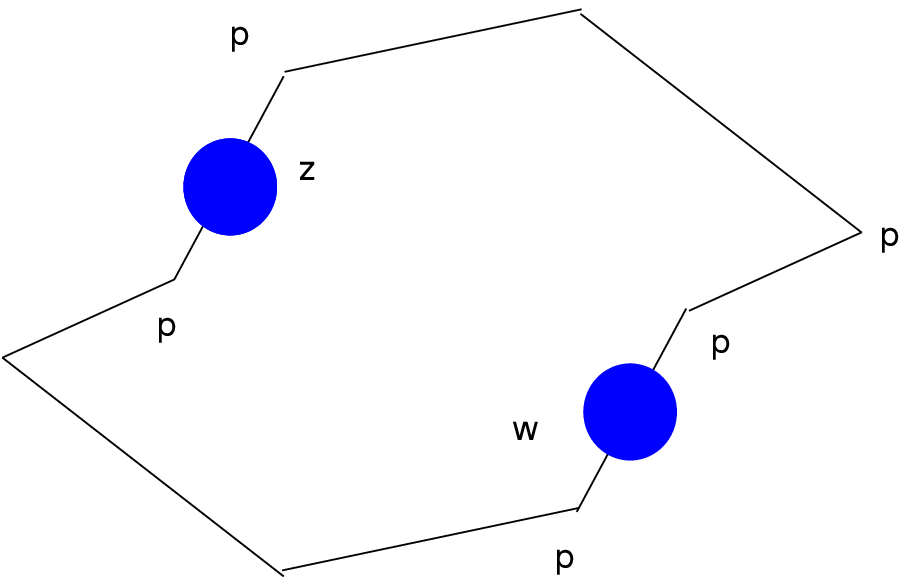}
    \caption{\label{fig:flat-p} Translation structure on $X\backslash p$}
    \end{figure}

For example, take an octagon $X$ with four pairs of parallel edges, see Figure~\ref{fig:flat}. 
Identifying the edges with the same labels by translation, $X$ becomes a
closed surface. All vertices are glued as one point $p$. By the topological Euler characteristic formula,        
the genus of $X$ is two. It is moreover a Riemann surface whose complex structure is induced from $\bbC$. 
Away from $p$ it admits an atlas of charts with transition functions given by translation: $z' = z + \text{constant}$, 
see Figure~\ref{fig:flat-p}.  

The differential $\omega = dz$ is well-defined and nowhere vanishing on $X\backslash p$, which further extends to the entire $X$. 
The angle at $p$ is $6\pi = 3 \cdot (2\pi)$, hence $\omega$ has a local expression $d (z^3) \sim z^2 dz$ at $p$. 
In summary, $\omega$ is an Abelian differential with a unique zero of order two on a Riemann surface of genus two. 

The above example illustrates the equivalence between translation surfaces and Abelian differentials in general. Given a translation surface, away from its saddle points, differentiating local coordinates provides a globally defined Abelian differential. Conversely, integrating an Abelian differential away from its zeros provides an atlas of charts whose transition functions are translation, because antiderivatives differ by constants. In addition, a saddle point $p$ has cone angle $(k+1)\cdot (2\pi)$ if and only if $\omega = d(z^{k+1})\sim z^k dz$ under a local coordinate $z$ at $p$, namely, if and only if $\omega$ has a zero of order $k$ at $p$. 

\begin{figure}[h]
    \centering
    \psfrag{a}{$a$}
    \psfrag{b}{$b$}
    \includegraphics[scale=0.25]{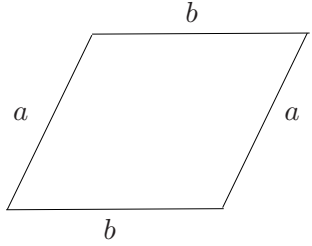}
    \caption{\label{fig:torus} A flat torus}
    \end{figure}

\begin{figure}[h]
    \centering
    \psfrag{a}{$a$}
    \psfrag{b}{$b$}
    \psfrag{c}{$c$}
    \psfrag{d}{$d$}
    \psfrag{e}{$e$}
    \includegraphics[scale=0.3]{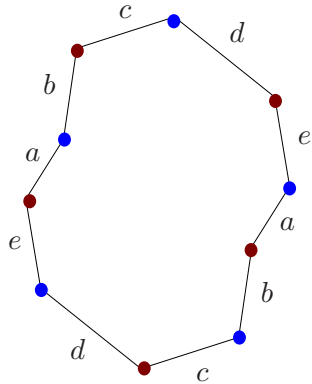}
    \caption{\label{fig:two-simple} A flat surface with two simple zeros}
    \end{figure}

Below we provide two more examples. Figure~\ref{fig:torus} represents a nowhere vanishing differential on a torus. Conversely every Abelian differential on a torus give rises to such a parallelogram presentation. Figure~\ref{fig:two-simple} represents an Abelian differential with two simple zeros on a Riemann surface of genus two.  

Note that Abelian differentials are sections of the \emph{canonical line bundle}, hence the study of translation surfaces is naturally connected to algebraic geometry. 

\subsection{Strata of Abelian differentials}
\label{subsec:strata}

We identify Riemann surfaces with smooth complex algebraic curves. Let $\MM_g$ be the \emph{moduli space of genus $g$ curves}. Let $\HH$ be the \emph{Hodge bundle} over $\MM_g$ whose fibers parameterize Abelian differentials on a fixed genus $g$ curve. 

Let $\mu = (m_1, \ldots, m_n)$ be a tuple of positive integers such that $\sum_{i=1}^n m_i = 2g-2$. We say that 
$\mu$ is a \emph{partition} of $2g-2$.\footnote{In Section~\ref{sec:higher} we will consider meromorphic differentials, where the entires $m_i$ are allowed to be negative.} 
Define a subset $\HH(\mu)$  of $\HH$ that parameterizes 
pairs $(X, \omega)$, where $X$ is a Riemann surface of genus $g$ and $\omega$ is an Abelian differential on $X$ such that the zero divisor of $\omega$ is of type 
$\mu$: 
$$ (\omega)_0 = m_1 p_1 + \cdots + m_n p_n. $$
We say that $\HH(\mu)$ is the \emph{stratum of Abelian differentials of type $\mu$}. 
Equivalently, $\HH(\mu)$ parameterizes translation surfaces with $n$ saddle points, each having cone angle 
$(m_i+1) \cdot (2\pi)$. The union of $\HH(\mu)$ over all partitions of $2g-2$ is the Hodge bundle $\HH$ (with the zero section removed). 

Take a basis $\gamma_1, \ldots, \gamma_{2g + n -1}$ of the relative homology $H_1(X, p_1, \ldots, p_n; \bbZ)$. Integrating 
$\omega$ over each $\gamma_i$ provides a local coordinate system for $\HH(\mu)$, called the \emph{period coordinates}. For instance, the complex vectors $v_1$, $v_2$, $v_3$, and $v_4$ in Figure~\ref{fig:flat} above are periods of a translation surface in $\HH(2)$. Under the polygon presentation, locally deforming the periods preserves the number of saddle points, their cone angles, and the way of edge identification. Consequently $\HH(\mu)$ is a $(2g+n-1)$-dimensional manifold.\footnote{More precisely it is an orbifold, because special translation surfaces can have extra automorphisms.} 

For special partitions $\mu$, $\HH(\mu)$ can be \emph{disconnected}. Kontsevich-Zorich 
(\cite{KontsevichZorich}) classified connected components of $\HH(\mu)$ for all $\mu$, where extra components arise due to hyperelliptic and spin structures. If a translation surface $(X, \omega)$ satisfies that $X$ is hyperelliptic, $(\omega)_0 = (2g-2)z$ or $(\omega)_0 = (g-1)(z_1+z_2)$, where 
$z$ is a Weierstrass point of $X$ in the former, or $z_1$ and $z_2$ are hyperelliptic conjugate in the latter, we say that $(X, \omega)$ is a \emph{hyperelliptic translation surface}.\footnote{Being a hyperelliptic translation surface not only requires $X$ to be hyperelliptic, but also imposes a condition on $\omega$.} For a nonhyperelliptic translation surface $(X, \omega)$, if 
$(\omega)_0 = 2 k_1 z_1 + \cdots + 2k_n z_n$, then the line bundle 
$\OO_X(\sum_{i=1}^n k_i z_i)$ 
is a square root of the canonical line bundle, namely, it is a \emph{theta characteristic}. Define its \emph{parity} by 
$h^0(X, \sum_{i=1}^n k_i z_i) \pmod{2}$, which is deformation invariant (\cite{Atiyah, Mumford}). 
A theta characteristic with its parity is called a \emph{spin structure}. In general, $\HH(\mu)$ can have \emph{up to three} connected components, distinguished by these hyperelliptic and spin structures. 

\subsection{$\GL_2^+(\bbR)$-action and affine invariant submanifolds}
\label{subsec:GL2}

\begin{figure}[h]
    \centering
    \psfrag{a}{$a$}
    \psfrag{b}{$b$}
    \psfrag{A}{$A$}
    \includegraphics[scale=0.4]{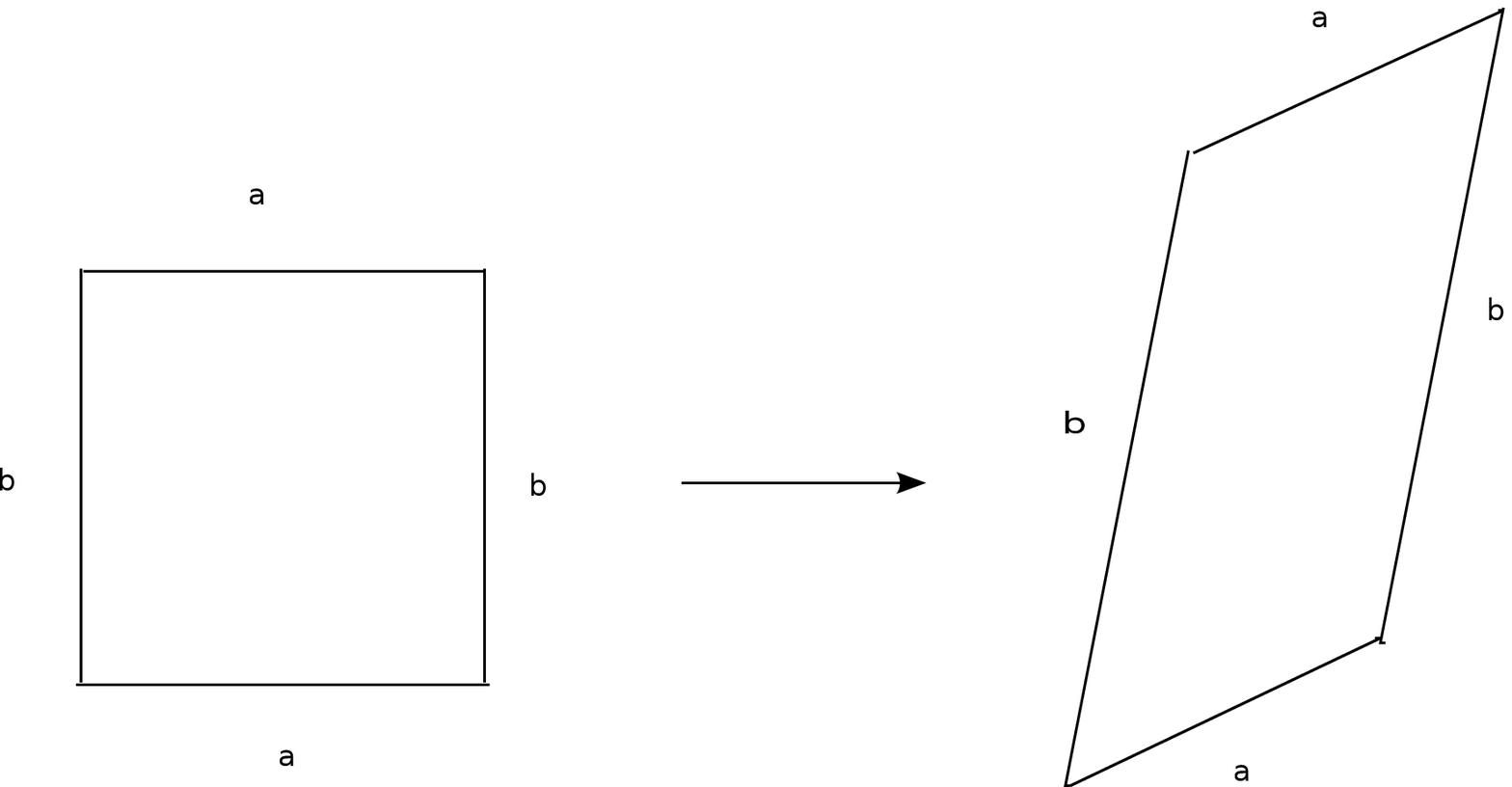}
    \caption{\label{fig:torus-GL} $\GL^{+}_2(\bbR)$-action on a flat torus}
\end{figure}

Given $(X, \omega) \in \HH$ and $A\in \GL^{+}_2(\bbR)$, varying the polygon presentation of $(X, \omega)$ by $A$ 
induces a $\GL^{+}_2(\bbR)$-action on $\HH$, which is called \emph{Teichm\"uller dynamics}. For example, the $\GL^{+}_2(\bbR)$-orbit of a flat torus consists of all parallelogram presentations, see Figure~\ref{fig:torus-GL}. 

The number and cone angles of saddle points are preserved under this action, hence the $\GL^{+}_2(\bbR)$-action descends to each stratum $\HH(\mu)$. Equip $\HH(\mu)$ the standard topology using its period coordinates. For almost all $(X, \omega) \in \HH(\mu)$, Masur and Veech (\cite{MasurInterval, VeechInterval}) showed that its $\GL^{+}_2(\bbR)$-orbit is \emph{equidistributed} in $\HH(\mu)$, hence the orbit closure is the whole stratum (or a connected component if the stratum is disconnected). For special $(X, \omega)$, however, its $\GL^{+}_2(\bbR)$-orbit closure can be a \emph{proper} subset of $\HH(\mu)$. Classifying 
$\GL^{+}_2(\bbR)$-orbit closures in $\HH(\mu)$ is a central question in Teichm\"uller dynamics. 

The recent breakthrough of Eskin-Mirzakhani-Mohammadi~(\cite{EskinMirzakhani, EskinMirzakhaniMohammadi}) showed that any $\GL^{+}_2(\bbR)$-orbit closure is an \emph{affine invariant submanifold} in $\HH(\mu)$, that is, locally it is a subspace of $\HH(\mu)$ cut out by real linear homogeneous equations of period coordinates.\footnote{Here the term ``affine'' is different from what it usually means in algebraic geometry. It refers to the linear structure on $\bbC^n$. Moreover, the closure of an orbit is taken under the standard topology in the Hodge bundle over the interior of the moduli space parameterizing smooth curves.} 
Filip (\cite{Filip}) further showed that all affine invariant submanifolds are algebraic varieties defined over $\overline{\bbQ}$. In particular, it means that affine invariant submanifolds can be defined and characterized purely in terms of algebraic conditions on the Jacobian. We will elaborate on these results in Section~\ref{sec:affine}. 

\subsection{Veech group}
\label{subsec:veech}

\begin{figure}[h]
    \centering
    \psfrag{a}{$a$}
    \psfrag{b}{$b$}
     \psfrag{c}{$c$}
    \psfrag{e}{$=$}
    \psfrag{A}{$\big(\begin{smallmatrix}
1 & 1 \\
0 & 1
\end{smallmatrix}\big)$}
    \includegraphics[scale=1.0]{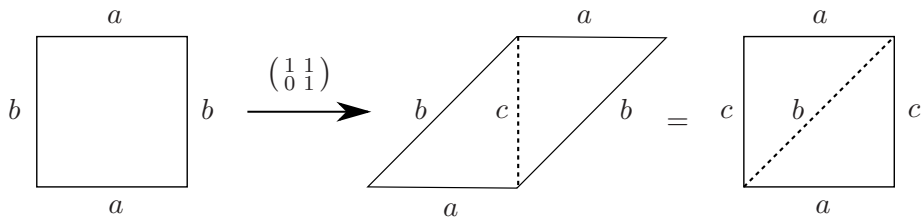}
    \caption{\label{fig:veech} An element of the Veech group}
\end{figure}

Let $(X,\omega) \in \HH(\mu)$ be a translation surface. Suppose a matrix $A\in \SL_2(\bbR)$ acts on $(X, \omega)$. 
If the resulting translation surface $A\cdot (X, \omega)$ is isomorphic to $(X,\omega)$, that is, if the polygon presentation 
of $A\cdot (X, \omega)$ can be cut into pieces and reassembled via translation to represent $(X,\omega)$, we say that 
$A$ is a \emph{stabilizer} of $(X, \omega)$. The subgroup of all stabilizers of $(X, \omega)$ in $\SL_2(\bbR)$ is called 
the \emph{Veech group}, and denoted by $\SL(X, \omega)$. For example, 
$A = \big(\begin{smallmatrix}
1 & 1 \\
0 & 1
\end{smallmatrix}\big)$ 
is in the Veech group of the square torus, see Figure~\ref{fig:veech}. 

A line segment under the flat metric that connects two zeros of $\omega$ (not necessarily distinct) is called a \emph{saddle connection}. Since the set of saddle connections of $(X, \omega)$ is preserved by a stabilizer, it follows that $\SL(X, \omega)$ is \emph{discrete}. Without loss of generality, suppose $(X, \omega)$ has a horizontal saddle connection. 
Let $g_t = \bigl(\begin{smallmatrix}
e^t & 0\\ 0 & e^{-t}
\end{smallmatrix} \bigr) 
 \in \SL_2(\bbR)$ act on $(X, \omega)$. As $t\to \infty$, the horizontal saddle connection becomes arbitrarily long, hence the 
$\SL_2(\bbR)$-orbit of $(X, \omega)$ is \emph{unbounded} in $\HH(\mu)$. Consequently $\SL(X, \omega)$ 
is \emph{not cocompact} in $\SL_2(\bbR)$. 

Adding the traces of all $A\in \SL(X, \omega)$ to $\bbQ$, we obtain a field extension of $\bbQ$, called the \emph{trace field} of 
$(X, \omega)$. The degree of the trace field over $\bbQ$ is bounded by the genus of $X$ 
(see \cite[Proposition 2.5]{MoellerSurvey}). 

\section{Teichm\"uller curves}
\label{sec:teich}

The Hodge bundle $\HH$ maps to the moduli space $\MM_g$ of smooth genus $g$ curves by forgetting the differentials: 
$(X, \omega) \mapsto X$. Note that the subgroup $\SO(2)$ acting on an Abelian differential amounts to rotating the corresponding flat surface, hence it does not change the underlying complex structure. Similarly scaling the size of a flat surface preserves the underlying complex structure. It follows that the projection of a $\GL^{+}_2(\bbR)$-orbit to $\MM_g$ factors through the upper half plane $\bbH$, and the induced map $\bbH \to \MM_g$ (or simply its image) is called a \emph{Teichm\"uller disk}. On rare occasions a Teichm\"uller disk forms an algebraic curve in $\MM_g$. In that case we call it a \emph{Teichm\"uller curve}.

\subsection{Properties of Teichm\"uller curves}
\label{subsec:teich}

Teichm\"uller curves are dimensionally minimal affine invariant submanifolds, which possess a number of fascinating properties. To name a few, a Teichm\"uller curve is a local \emph{isometry} from a curve to $\MM_g$ under the 
Kobayashi/Teichm\"uller metric (\cite{SmillieWeiss, VeechTeich}). The union of all Teichm\"uller curves is \emph{dense} in moduli spaces (\cite{EskinOkounkov, ChenRigid}). McMullen (\cite{McMullenRigid}) proved that Teichm\"uller curves are \emph{rigid}, hence they are defined over number fields.\footnote{Here the rigidity means as a map it does not deform.} Conversely, Ellenberg-McReynolds (\cite{EllenbergMcReynolds}) showed that every curve over a number field is birational to a Teichm\"uller curve over $\bbC$. If $(X, \omega)$ generates an algebraically primitive Teichm\"uller curve (see Section~\ref{subsec:teich-class}), M\"oller (\cite{MoellerTorsion})
showed that the difference of any two zeros of $\omega$ is a torsion in the Jacobian of $X$. M\"oller (\cite{MoellerHodge})
also analyzed the variation of Hodge structures associated to a Teichm\"uller curve and deduced that it parameterizes curves whose Jacobians have \emph{real multiplication}. Teichm\"uller curves are never complete in $\MM_g$, but the closure of a Teichm\"uller curve only intersects certain boundary divisor of the Deligne-Mumford compactification $\BM_g$ (see Section~\ref{subsec:slope}).\footnote{Here we restrict to Teichm\"uller curves generated by Abelian differentials. Special quadratic differentials can also generate Teichm\"uller curves, see Section~\ref{subsec:quad}, but they may intersect other boundary divisors of $\BM_g$.} 

\subsection{Square-tiled surfaces}
\label{subsec:square-tiled}

We first show that Teichm\"uller curves exist. One type of Teichm\"uller curves arises from certain \emph{branched covering construction}. Let $\mu = (m_1, \ldots, m_n)$ be a partition of $2g-2$. Consider a branched cover $\pi: X \to E$ such that  
\begin{itemize}
\item $\deg \pi = d$, 

\item $g(X) = g$, 

\item $E$ is the square torus, 

\item $\pi$ has a \emph{unique} branch point $q\in E$,

\item $\pi$ has $n$ ramification points $p_1, \ldots, p_n$ over $q$, each with ramification order $m_i$.  
\end{itemize}

\begin{figure}[h]
    \centering
    \psfrag{a}{$a$}
    \psfrag{b}{$b$}
    \psfrag{C}{$X$}
    \psfrag{c}{$c$}
    \psfrag{d}{$d$}
    \psfrag{E}{$E$}
    \psfrag{u}{$u$}
    \psfrag{v}{$v$}
     \psfrag{e}{$5:1$}
    \includegraphics[scale=0.5]{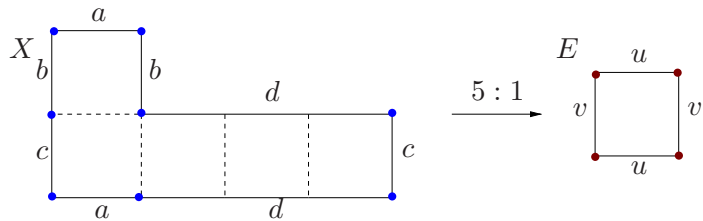}
    \caption{\label{fig:cover} A square-tiled surface in $\HH(2)$}
    \end{figure}

Let $\omega = \pi^{*} (dz)$, where $z$ is the standard coordinate on $E$. Then by the Riemann-Hurwitz formula 
$$(\omega)_0 = m_1 p_1 + \cdots + m_n p_n,$$ 
hence $(X, \omega) \in \HH(\mu)$. Such $(X, \omega)$ are called \emph{square-tiled surfaces} (or \emph{origami}). For example, Figure~\ref{fig:cover} exhibits a degree $5$ and genus $2$ branched cover of $E$ with a unique ramification point of order $2$. The resulting $(X, \omega)$ belongs to the stratum $\HH(2)$. 

\begin{figure}[h]
    \centering
      \includegraphics[scale=0.5]{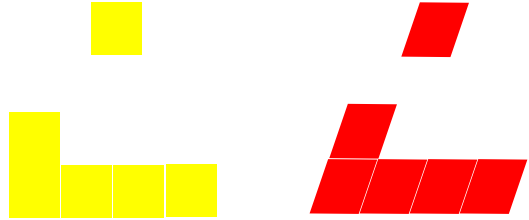}
    \caption{\label{fig:cover-GL} $\GL^{+}_2(\bbR)$-action on a square-tiled surface}
    \end{figure}

The $\GL^{+}_2(\bbR)$-action on a square-tiled surface 
amounts to varying the shape of the square, which is exchangeable with varying the shape of the flat torus first, see Figure~\ref{fig:cover-GL} for an example. Therefore, the Teichm\"uller disk generated by a square-tiled surface corresponds to the one-dimensional \emph{Hurwitz space} parameterizing degree $d$, genus $g$ connected covers of all elliptic curves with a unique branch point of ramification type $\mu$. Since Hurwitz spaces are algebraic varieties, it follows that the $\GL^{+}_2(\bbR)$-orbit of a square-tiled surface gives rise to a Teichm\"uller curve.\footnote{When $d$, $g$, and $\mu$ are fixed, the Hurwitz space can still be disconnected. In that case each of its connected components gives a Teichm\"uller curve.} Since $d$ can be arbitrarily large, this way we indeed obtain \emph{infinitely many} Teichm\"uller curves in $\HH(\mu)$.  
        
\begin{figure}[h]
    \centering
      \includegraphics[scale=1.0]{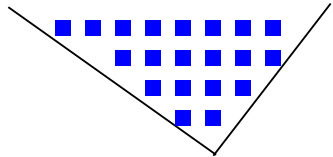}
    \caption{\label{fig:lattice-square} Square-tiled surfaces as lattice points in a stratum}
    \end{figure}

Square-tiled surfaces correspond to \emph{lattice points} under period coordinates of a stratum (\cite[Lemma 3.1]{EskinOkounkov}), see Figure~\ref{fig:lattice-square}. 
If $\pi: X\to E$ is a square-tiled surface of type $\mu = (m_1, \ldots, m_n)$, for any $\gamma\in H_1(X, p_1, \ldots, p_n; \bbZ)$, $\pi(\gamma)$ represents a closed loop in $E$, because 
$p_1, \ldots, p_n$ all map to the unique branch point. Therefore, 
$$\int_{\gamma}\pi^{*}(dz) = \int_{\pi_{*}\gamma} dz \in \bbZ \oplus \bbZ[i]. $$
Conversely if all relative periods of $(X,\omega)$ are lattice points, the map $\pi: X\to \bbC/\bbZ \oplus \bbZ[i]$  
induced by 
$$x\mapsto \int_{b}^x \omega$$ 
is well-defined, where $b$ is a fixed base point, hence realizing $(X,\omega)$ as a square-tiled surface. This explains density of 
the union of Teichm\"uller curves in moduli spaces. It also provides an approach for analyzing \emph{volume growths} of the strata of Abelian differentials by counting the number of such square-tiled surfaces (\cite{EskinOkounkov}).  

Gutkin-Judge (\cite{GutkinJudge}) showed that the Veech group of a translation surface is 
{\em commensurable with $\SL(2,\bbZ)$} if and only if the surface is tiled by parallelograms. In this sense Teichm\"uller curves generated by square-tiled surfaces are called \emph{arithmetic} Teichm\"uller curves. 
However, there exist Teichm\"uller curves of other type that are \emph{not} generated by square-tiled surfaces. In the next section we will survey known results on the classification of Teichm\"uller curves. 

\subsection{Classification of Teichm\"uller curves}
\label{subsec:teich-class}

By definition, the moduli space $\MM_{1,1}$ of elliptic curves is a very first example of Teichm\"uller curves. Arithmetic Teichm\"uller curves generated by square-tiled surfaces are coverings of $\MM_{1,1}$. On the contrary, a Teichm\"uller curve is called (geometrically) \emph{primitive} if it does not arise from a curve in a moduli space of lower genus via such a covering construction.

It is a challenging task to find examples of primitive Teichm\"uller curves. In genus two, Calta and McMullen (\cite{Calta, McMullenWeier})
independently classified primitive Teichm\"uller curves in $\HH(2)$. They found infinitely many primitive Teichm\"uller curves whose constructions have several incarnations. In terms of flat geometry, a translation surface $(X, \omega)$ generating a primitive Teichm\"uller curve in $\HH(2)$ possesses an $L$-shaped polygon presentation whose edges satisfy relations in a real quadratic field. In terms of algebraic geometry, the Jacobian of $X$ admits real multiplication, and the corresponding Teichm\"uller curve lies on a Hilbert modular surface. However for the other stratum $\HH(1,1)$ in genus two, McMullen (\cite{McMullenTorsion}) proved that it contains a unique primitive Teichm\"uller curve generated by the regular decagon with parallel edges identified. McMullen (\cite{McMullenPrym}) further generalized the method of real multiplication by using Prym varieties and discovered infinitely many primitive Teichm\"uller curves in genus three and four. Using quotients of abelian covers of $\bbP^1$ by a finite group, in each genus Bouw-M\"oller (\cite{BouwMoeller})
constructed (finitely many) primitive Teichm\"uller curves, generalizing earlier constructions of Veech and Ward (\cite{VeechTriangle, Ward}).

Besides the above results, to date only few sporadic examples of primitive Teichm\"uller curves are known. It is natural to ask if there are only finitely many primitive Teichm\"uller curves in a given stratum $\HH(\mu)$. There is a stronger notion of primitivity called \emph{algebraically primitive}, if the trace field of a translation surface $X$ in the Teichm\"uller curve has degree equal to the genus of $X$ (see Section~\ref{subsec:veech}). For example, the trace field of a square-tiled surface is $\bbQ$, hence it is far from being algebraically primitive. Indeed, algebraically primitive Teichm\"uller curves are geometrically primitive, but the converse is not always true (see \cite[Section 5.1]{MoellerSurvey}). 

Finiteness results of algebrically primitive Teichm\"uller curves have been established in various cases. M\"oller (\cite{MoellerFinite})
proved that the hyperelliptic component of $\HH(g-1, g-1)$ contains finitely many algebrically primitive Teichm\"uller curves. The strategy is to track the degeneration of flat surfaces along an algebrically primitive Teichm\"uller curve in the horizontal and vertical directions and use it to bound the torsion order of the difference of the two zeros. By studying the boundary of the locus of curves with real multiplication, Bainbridge-M\"oller (\cite{BainbridgeMoellerReal})
showed finiteness of algebrically primitive Teichm\"uller curves in $\HH(3,1)$. Bainbridge-Habegger-M\"oller (\cite{BainbridgeHabeggerMoeller}) further established finiteness of algebrically primitive Teichm\"uller curves for all strata in genus three by a mix of techniques, including the Harder-Narasimhan filtration of the Hodge bundle over Teichm\"uller curves and height bounds for the boundary points of Teichm\"uller curves. Matheus-Wright (\cite{MatheusWright})
proved finiteness of algebrically primitive Teichm\"uller curves in the minimal stratum $\HH(2g-2)$ for each prime genus $g \geq 3$. Their approach is to study orthogonality of Hodge-Teichm\"uller (real) planes in the Hodge bundle that respect the Hodge decomposition along a Teichm\"uller curve. Matheus-Nguyen-Wright (\cite{NguyenWright})
 further showed that there are at most finitely many non-arithmetic Teichm\"uller curves in the hyperelliptic component of $\HH(4)$. 

One can also study Teichm\"uller curves contained in an affine invariant submanifold of a stratum. 
Lanneau-Nguyen (\cite{LanneauNguyen}) showed that there are at most finitely many closed $\GL^{+}_2(\bbR)$-orbits (including primitive Teichm\"uller curves) in certain Prym loci in genus three. Lanneau-Nguyen-Wright (\cite{LanneauNguyenWright})
further proved finiteness of closed $\GL^{+}_2(\bbR)$-orbits in each non-arithmetic rank $1$ affine invariant submanifold (see Section~\ref{subsec:affine-class}). Apisa (\cite{Apisa}) 
showed finiteness of algebrically primitive Teichm\"uller curves in the hyperelliptic components of each minimal stratum in $g > 2$, as a byproduct of his classification of affine invariant submanifolds in the hyperelliptic components (see Section~\ref{subsec:affine-class}). 

Although a complete classification of Teichm\"uller curves is still missing, the seminal work of Eskin-Mirzakhani-Mohammadi~(\cite{EskinMirzakhani, EskinMirzakhaniMohammadi}) on the structure of $\GL^{+}_2(\bbR)$-orbit closures provides us a powerful tool. Indeed some of the above results are built on their work. We hope the classification problem of Teichm\"uller curves (and in general affine invariant submanifolds) can be resolved in the next few decades. 

\subsection{Slope, Siegel-Veech constant, and Lyapunov exponent}
\label{subsec:slope}

In this section we discuss several invariants of Teichm\"uller curves and describe a relation between them. 

\begin{figure}[h]
    \centering
    \psfrag{L}{$L$}
      \includegraphics[scale=0.8]{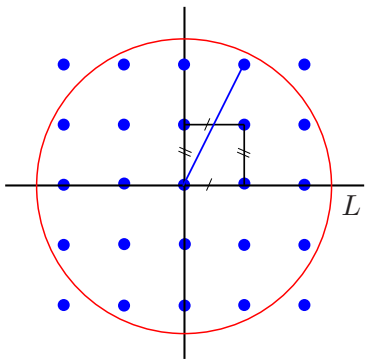}
    \caption{\label{fig:lattice-saddle} Saddle connections on the standard torus}
    \end{figure}

First, the behavior of geodesics on translation surfaces is related to \emph{billiards in polygons}. One can study various counting problems from this viewpoint. Recall that a saddle connection is a line segment connecting two saddle points. Consider counting saddle connections with bounded lengths. For instance, consider the standard torus formed by identifying parallel edges of the unit square and marked at the origin. The number of saddle connections of length $< L$ (counting with direction) equals the number of lattice points in the disk of radius $L$, see Figure~\ref{fig:lattice-saddle}. 
It has asymptotically \emph{quadratic growth} $\sim \pi L^2$. The leading term $\pi$ is an example of so-called \emph{Siegel-Veech constant}. In general, Siegel-Veech constants relate quadratic growth rates of finite-length trajectories on a translation surface to the volume of the corresponding $\SL_2(\RR)$-orbit. In what follows we will concentrate on one type of Siegel-Veech constants, called the \emph{area} Siegel-Veech constant, which has a connection to dynamics as well as intersection theory on moduli space. 

\begin{figure}[h]
    \centering
    \psfrag{h1=1}{$h_1 = 1$}
    \psfrag{l1=1}{$w_1 = 1$}
    \psfrag{h2=1}{$h_2 = 1$}
    \psfrag{l2=4}{$w_2 = 4$}
     \includegraphics[scale=0.5]{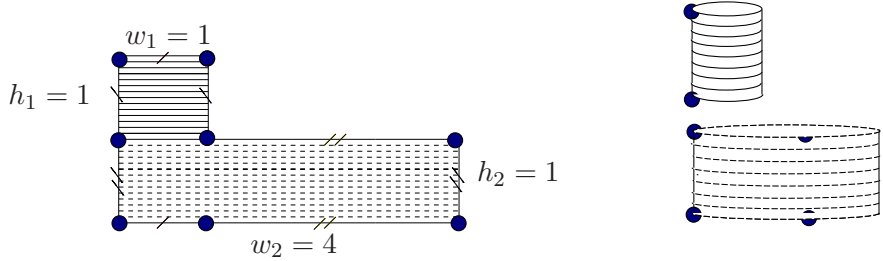}
    \caption{\label{fig:cylinder} Two horizontal cylinders on a square-tiled surface}
    \end{figure}

A real geodesic on a translation surface $(X, \omega)$ is called \emph{regular}, if it does not pass through any saddle point, namely, it does not contain any zero of $\omega$. Vary a closed regular geodesic in parallel until it hits a saddle point on both ends. The union of those geodesics fill a \emph{cylinder} $\cyl$, whose boundary circles contain saddle points. For example, the square-tiled surface in Figure~\ref{fig:cylinder} has two cylinders in the horizontal direction. 
Let $w$ and $h$ be the \emph{width} and  \emph{height} of a cylinder $\cyl$, respectively. The \emph{area} of the cylinder is 
$\Area(\cyl) = w \cdot h$. 

For $L > 0$, define 
$$N(X, L) = \frac{1}{\Area(X)} \cdot \sum_{\cyl\subset X\atop w(\cyl) < L} \Area(\cyl).$$ 
Veech and Eskin-Masur (\cite{VeechSiegel, EskinMasur}) showed that for any $\GL^{+}_2(\bbR)$-orbit closure $\calN$, there exists a constant $c$ such that 
$$ \frac{\pi}{3}\cdot \lim_{L \to \infty} \frac{N(X,L)}{L^2} = c $$
for \emph{generic} $(X, \omega)\in \calN$. The constant 
$c$ is called the \emph{area Siegel-Veech constant} of $\calN$.\footnote{Here our definition of the Siegel-Veech constant differs from its usual definition by a scalar multiple of $\pi^2/3$. We choose to normalize it this way for the convenience of describing its relation with slope and Lyapunov exponent.} Its value depends on $\calN$. 

If $\calN$ is the orbit of a square-tiled surface, namely, if its projection to $\MM_g$ is an arithmetic Teichm\"uller curve, there is a \emph{combinatorial} way to calculate its area Siegel-Veech constant. Let $N$ be the number of square-tiled surfaces in $\calN$. In other words, $N$ is the \emph{Hurwitz number} that counts branched covers of a fixed torus with a unique branch point and prescribed ramification type. Take all square-tiled surfaces in $\calN$. For each one, consider all of its horizontal cylinders. Sum up $h/w$ (the \emph{modulus} of a cylinder) over all of them. Denote the total sum by $M$. For example, the square-tiled surface in Figure~\ref{fig:cylinder} above contributes $ \frac{h_1}{w_1} + \frac{h_2}{w_2} = \frac{1}{1} + \frac{1}{4}$
to the sum $M$. 

For an arithmetic Teichm\"uller curve generated by a square-tiled surface,  
$$ c = \frac{M}{N}, $$
see~\cite[Appendix B]{EKZ}. 
The idea behind the formula is that $c$ measures the average number of cylinders weighted by their moduli $h/w$ in 
$\calN$, where 
$$\frac{h}{w} = \frac{hw}{w^2} = \frac{\Area(\cyl)}{w^2}, $$
hence it is related to the quadratic growth rate of $N(X, L)$ in this case. 

For low degree $d$ and low genus $g$, using monodromy of branched covers one can calculate $N$ and $M$ explicitly. Nevertheless, the enumeration of $N$ 
as $d$ and $g$ increase is a highly non-trivial problem in symmetric group representations. Eskin-Okounkov (\cite{EskinOkounkov}) 
analyzed the asymptotic behavior of $N$ and calculated the \emph{volume} growth of strata of Abelian differentials.  
The enumeration of $M$ is more complicated for large $d$ and $g$. Joint with M\"oller and Zagier (\cite{ChenMoellerZagier}) we are able to understand the asymptotic growth of $M$ and hence $c$ for arithmetic Teichm\"uller curves by using techniques of shifted symmetric functions and quasimodular forms. 

\begin{figure}[h]
    \centering
    \psfrag{i}{$i$}
    \psfrag{a}{$\Delta_i$:}
    \psfrag{j}{$g-i$}
     \psfrag{k}{$g-1$}
     \psfrag{1}{$i-1$}
    \psfrag{b}{$\Delta_0$:}
    \psfrag{c}{$\Delta_0\cap \Delta_i$:}
      \includegraphics[scale=0.3]{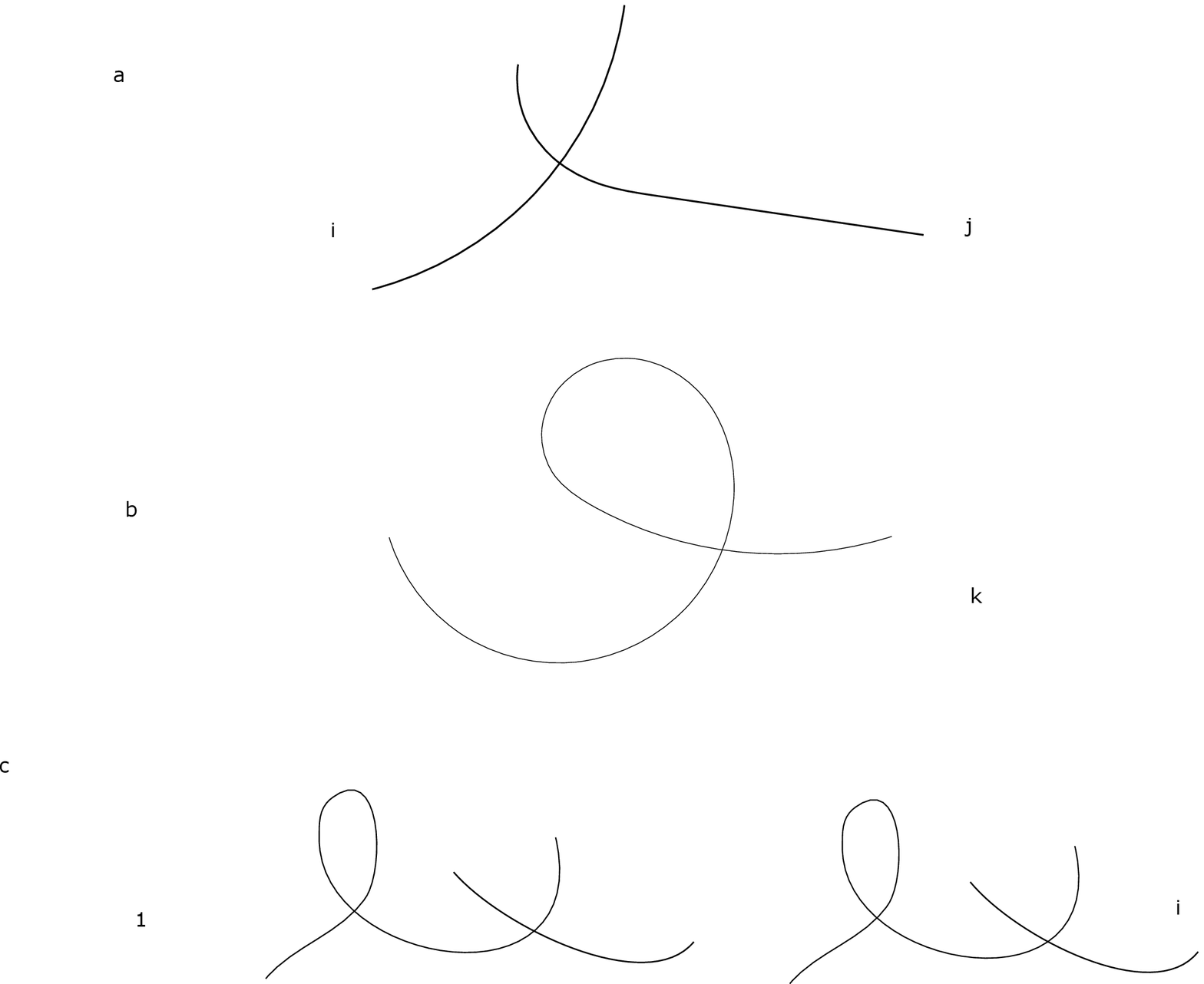}
    \caption{\label{fig:boundary} Curves in the boundary of $\BM_g$}
    \end{figure}

Next we define an important index associated to a one-dimensional family of stable curves. The \emph{boundary} $\Delta$ of the 
Deligne-Mumford moduli space $\BM_g$ parameterizes \emph{stable nodal curves}, where 
$\Delta = \bigcup\limits_{i=0}^{[g/2]} \Delta_i$ is a union of irreducible \emph{boundary divisors} $\Delta_i$. 
General points of $\Delta_i$ parameterize nodal curves of a given topological type, and they can further degenerate, 
see Figure~\ref{fig:boundary}. 

\begin{figure}[h]
    \centering
    \psfrag{B}{$B$}
      \includegraphics[scale=0.2]{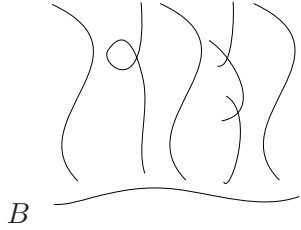}
    \caption{\label{fig:family} A one-dimensional family of stable curves}
    \end{figure}

The Hodge bundle $\HH$ over $\MM_g$ extends to a rank $g$ bundle $\BHH$ over $\BM_g$, whose fiber over a stable curve $X$ parameterizes sections of the dualizing line bundle $K_X$. Geometrically speaking, the fiber of $\BHH$ over $X$ can be identified with the space of \emph{stable differentials} that have at worst simple pole at each node of $X$ with opposite residues on the two branches of a node (\cite[Chapter 3.A]{HarrisMorrison}). Denote by $\lambda$ the first Chern class of $\BHH$ over $\BM_g$. Given a one-dimensional family $B$ of stable genus $g$ curves, define its \emph{slope} by  
$$ s = \frac{\deg \Delta|_B}{\deg \lambda|_B}. $$
Morally speaking, $\deg \Delta|_B$ counts the number of nodes (with multiplicity) and $\deg \lambda|_B$ measures the variation of complex structures, see Figure~\ref{fig:family}.  

Let $\mu = (m_1, \ldots, m_n)$ be a partition of $2g-2$. Define 
$$\kappa_{\mu} = \frac{1}{12}\sum_{i=1}^{n} \frac{m_i (m_i+2)}{m_i+1},$$ 
which depends on $\mu$ only. For a Teichm\"uller curve in $\HH(\mu)$, despite that its area Siegel-Veech constant and slope are defined in different contexts, they 
determine each other (\cite{ChenThesis, ChenRigid, ChenMoellerAbelian}): 
$$ s = \frac{12c}{c + \kappa_{\mu}}. $$ 
    
\begin{figure}[h]
    \centering
    \psfrag{h}{$h$}
    \psfrag{w}{$w$}
    \psfrag{n}{a node}
 \includegraphics[scale=0.25]{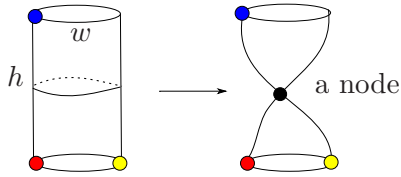}
    \caption{\label{fig:delta} Shrinking the core curve of a cylinder}
    \end{figure}

The upshot to prove the formula consists of the following. First, the constant $\kappa_\mu$ corresponds to the Miller-Morita-Mumford \emph{$\kappa$-class}. Next, a cylinder with modulus $h/w$ contributes $h/w$ to the intersection of the Teichm\"uller curve with the boundary $\Delta$, see Figure~\ref{fig:delta}. Moreover, the numerical relation $\displaystyle 12\lambda \equiv \kappa+ \Delta$ holds on $\BM_g$.  

Finally we introduce an important dynamical invariant for an affine invariant submanifold. The diagonal subgroup $\bigl(\begin{smallmatrix}
e^t & 0\\ 0 & e^{-t}
\end{smallmatrix} \bigr)$ 
defines the \emph{Teichm\"uller geodesic flow} on the Hodge bundle $\HH$ over an affine invariant submanifold $\calN$.  
Since $\HH$ is of rank $g$, there are $g$ nonnegative \emph{Lyapunov exponents} 
$\lambda_1 \geq \cdots \geq \lambda_g \geq 0$ 
as logarithms of mean eigenvalues of monodromy of $\HH$ along the flow on $\calN$. Morally speaking, Lyapunov exponents  
measure the separation rates of infinitesimally closed trajectories (see \cite{Zorich} for more details). Denote the \emph{sum} of the Lyapunov exponents by 
$$L = \lambda_1 + \cdots + \lambda_g.$$ 

In a seminal work, Eskin-Kontsevich-Zorich (\cite{EKZ})
proved that the sum of Lyapunov exponents and the area Siegel-Veech constant determine each other for \emph{any} affine invariant submanifold: 
$$ L = c + \kappa_{\mu}. $$
As a corollary, we thus obtain that for a Teichm\"uller curve, 
any one of the three numbers $s$, $c$, and $L$ determine the other two. 

As an application, one can deduce a non-varying phenomenon of Teichm\"uller curves in low genus. By computer experiments, Kontsevich and Zorich observed that for many strata $\HH(\mu)$ in low genus, \emph{all} Teichm\"uller curves in the same stratum (component) have \emph{non-varying} sums of Lyapunov exponents. They came up with a conjectural list of such strata: 
\begin{itemize}
\item $g =2$: $\HH(1,1)$, $\HH(2)$. 

\item $g=3$: $\HH(4)$, $\HH(3,1)$, $\HH(2,2)$, $\HH(2,1,1)$.

\item $g=4$: $\HH(6)$, $\HH(5,1)$, $\HH(4,2)$, $\HH(3,3)$, $\HH(3,2,1)$, $\HH^{\odd}(2,2,2)$.

\item $g=5$:  $\HH(8)$, $\HH^{\odd}(6,2)$, $\HH(5,3)$, $\HH^{\hyp}(4,4)$. 
\end{itemize}

This phenomenon disappears when $g\geq 6$ except for the \emph{hyperelliptic} strata. 

Joint with M\"oller we proved Kontsevich-Zorich's conjecture (\cite{ChenMoellerAbelian}). Let us use $\HH(3,1)$ as an example to explain our method. 
For $(X, \omega)\in \HH(3,1)$, it is easy to observe that $X$ is \emph{non-hyperelliptic}. The same holds for any nodal curve $X$ contained in the \emph{boundary} of a Teichm\"uller curve $\TT$ in $\HH(3,1)$. Let 
$\Hyp$ be the closure of the locus of hyperelliptic curves in $\BM_3$, which is a divisor with class 
$$\Hyp \equiv 9 \lambda - \Delta_0 - 3\Delta_1, $$
see \cite[Chapter 3.H]{HarrisMorrison}. Since $\TT$ is disjoint with $\Hyp$, we conclude that $\TT \cdot \Hyp = 0$. Moreover, 
$\TT\cdot \Delta_i = 0$ for $i>0$, because shrinking the core curve of a cylinder yields a node of type $\Delta_0$ only, see 
Figure~\ref{fig:delta} above. By $\TT \cdot (9 \lambda - \Delta) = 0$, we see that the slope   
$$ s = \displaystyle\frac{\deg \Delta|_\TT}{\deg \lambda|_\TT} = 9.$$ 
The above calculation is independent of the choice of Teichm\"uller curves in $\HH(3,1)$. In summary, all Teichm\"uller curves in $\HH(3,1)$ have slope $9$. By the slope, Siegel-Veech, and Lyapunov exponent formula, sums of Lyapunov exponents are the same for all Teichm\"uller curves in $\HH(3,1)$. 

In general, the strategy is to find a geometrically meaningful \emph{divisor} $D$ in $\BM_g$ (or in $\BM_{g,n}$ by marking the zeros of $\omega$) such that $D$ is \emph{disjoint} with all Teichm\"uller curves in $\HH(\mu)$. The divisor class of $D$ determines the invariants of those Teichm\"uller curves. 

We remark that Yu-Zuo (\cite{YuZuo}) gave another novel proof of Kontsevich-Zorich conjecture using the \emph{Harder-Narasimhan filtration} of the Hodge bundle over Teichm\"uller curves. Their idea is the following. Suppose $f: \calX\to \TT$ is the universal curve over a Teichm\"uller curve. Let $S_1, \ldots, S_n$ be the disjoint sections of zeros of $(X, \omega)\in \TT$. For a tuple of integers $a_1, \ldots, a_n$, if $h^0(X, \sum_{i=1}^n a_i z_i) = k$ for all $X$, then the direct image sheaf 
$f_{*}\OO_{\calX}(\sum_{i=1}^n a_i S_i)$ is a vector bundle of rank $k$ on $\TT$. 

Use $\HH(3,1)$ as an example again. Since 
$3z_1 + z_2 \sim K_X$, one checks that $h^0(X, z_1 + z_2 ) = 1$ and $h^0(X, 2z_1 + z_2) = 2$ for all $X \in \TT$ (including the boundary points). Hence we obtain a filtration 
$$ f_{*}\OO_{\calX}(S_1 + S_2)\subset f_{*}\OO_{\calX}(2S_1 + S_2)\subset  f_{*}\OO_{\calX}(3S_1 + S_2), $$
where the last term can be identified with the Hodge bundle twisted by the generating differential of $\TT$. Then the sum $L$ of 
Lyapunov exponents of $\TT$ follows from Chern class calculations of the subbundles in the filtration. 

We also remark that the slopes of Teichm\"uller curves can provide information to understand the \emph{cone of effective divisors} 
of $\BM_g$. The idea is that any effective divisor on $\BM_g$ cannot contain all arithmetic Teichm\"uller curves, because their union is \emph{dense} in $\BM_g$. Hence the limit of slopes of arithmetic Teichm\"uller curves as the degree of the coverings approach infinity bounds the numerical class of the divisor. We refer to the survey (\cite{ChenFaraksMorrison}) for more details on the effective cone of moduli spaces. 

\subsection{Miscellaneous}
\label{subsec:misc}

In this section we collect various results with a flavor in algebraic geometry that are related to the preceding discussions. 

McMullen (\cite{McMullenSpin}) classified connected components of arithmetic Teichm\"uller curves in $\HH(2)$. The result is that an arithmetic Teichm\"uller curve in $\HH(2)$ is either connected or has two components. The special case when the degree of the coverings is prime was previously established by 
Hubert-Leli\`{e}ver (\cite{HubertLelievre}). In general, classifying connected components of arithmetic Teichm\"uller curves is a wide open question. From the viewpoint of algebraic geometry, it amounts to classifying connected components of the Hurwitz space of torus coverings with a unique branch point and prescribed ramification type. Bainbridge (\cite{BainbridgeEuler})
calculated Euler characteristics of Teichm\"uller curves in $\HH(2)$. The idea is to calculate the fundamental classes of these Teichm\"uller curves in certain compactifications of Hilbert modular surfaces. Mukamel (\cite{Mukamel})
further determined the number and type of orbifold points of these Teichm\"uller curves. Kumar-Mukamel (\cite{KumarMukamel})
described algebraic models of Teichm\"uller curves in genus two. Weiss (\cite{Weiss})
calculated the volume of certain twisted Teichm\"uller curves on Hilbert modular surfaces and partially classified their connected components.  

Siegel-Veech constants and Lyapunov exponents can be defined similarly for any affine invariant submanifold, say the strata themselves. Eskin-Masur-Zorich (\cite{EskinMasurZorich}) analyzed the principal boundary of the moduli space of Abelian differentials and gave a recursive method to compute Siegel-Veech constants of the strata for any Siegel-Veech configuration. For the cylinder configuration, the limit of area Siegel-Veech constants of arithmetic Teichm\"uller curves in a stratum as the degree of the coverings approach infinity equals the area Siegel-Veech constant of the stratum (\cite[Appendix A]{ChenRigid}). 

Previously we discussed the sum of Lyapunov exponents. For arithmetic Teichm\"uller curves generated by cyclic covers of $\bbP^1$ branched at four points, Eskin-Kontsevich-Zorich (\cite{EKZ2}) calculated all individual Lyapunov exponents. Yu (\cite{Yu})
conjectured that the polygon of Lyapunov spectrum bounds the Harder-Narasimhan polygon on Teichm\"uller curves, and a proof of this conjecture is recently announced by Eskin-Kontsevich-M\"oller-Zorich. The fundamental work of Forni (\cite{Forni})
showed that for a stratum of Abelian differentials no Lyapunov exponent vanishes. The fundamental work of Avila and Viana (\cite{AvilaViana}) further showed that for a stratum the Lyapunov spectrum is simple, that is, the strict inequality 
$\lambda_i > \lambda_{i+1}$ holds for all $i$. For arithmetic Teichm\"uller curves generated by square-tiled surfaces,  
Matheus-M\"oller-Yoccoz (\cite{MatheusMoellerYoccoz}) gave a Galois-theoretical criterion for the simplicity of their Lyapunov spectra.  

The Lyapunov spectrum can degenerate for a special affine invariant submanifold. If $\lambda_2 = \cdots = \lambda_g = 0$, we say that the Lyapunov spectrum is completely degenerate. 
Forni and Forni-Matheus-Zorich (\cite{ForniDegenerate, ForniMatheusZorich})
found examples of Teichm\"uller curves with completely degenerate Lyapunov spectrum in genus three and four, respectively. By studying Teichm\"uller curves that 
are also Shimura curves, M\"oller (\cite{MoellerST}) showed that those examples are the only Teichm\"uller curves with completely degenerate Lyapunov spectrum with possible exceptions in genus five. Aulicino (\cite{Aulicino})
showed that an affine invariant submanifold with completely degenerate Lyapunov spectrum can only be an arithmetic Teichm\"uller curve in genus at most five, and he further established finiteness of Teichm\"uller curves with completely degenerate Lyapunov spectrum. Filip (\cite{FilipDegenerate})
described all situations when an affine invariant submanifold can possess a zero Lyapunov exponent by analyzing monodromy of the corresponding Kontsevich-Zorich cocycle. As a higher dimensional analogue of families of curves, Filip (\cite{FilipK3})
considered families of K3 surfaces whose second cohomology groups form a local system, and showed that their top Lyapunov exponents are always rational. 

To conclude this section we mention two problems that have broad connections to other fields. Kontsevich-Zorich (\cite{KontsevichZorich2}) conjectured that every connected component of the strata is $K(\pi, 1)$, that is, it has a contractible universal cover and its fundamental group is commensurable with certain mapping class group. For all strata in genus three except $\HH(1^4)$, Looijenga-Mondello (\cite{LooijengaMondello}) determined their (orbifold) fundamental groups by analyzing geometry of canonical curves. Kontsevich-Soibelman (\cite{KontsevichSoibelman}) speculated that moduli spaces of differentials can be identified with
moduli spaces of certain stability conditions, where the $\GL_2^{+}(\bbR)$-action and saddle connections are analogues 
of the central charge in the theory of stability conditions. The seminal work of Bridgeland-Smith (\cite{BridgelandSmith}) established this identification for moduli spaces of quadratic differentials with simple zeros by relating the finite-length trajectories of such quadratic differentials to the stable objects of the corresponding stability condition. 

\section{Affine invariant submanifolds}
\label{sec:affine}

In the preceding section we discussed Teichm\"uller curves as examples of affine invariant submanifolds. 
In this section we consider affine invariant submanifolds in general. 

\subsection{Structure of affine invariant submanifolds}
\label{subsec:structure}

The celebrated Ratner's orbit closure theorem in ergodic theory says that the closures of orbits of unipotent flows on the quotient of a Lie group by a lattice are homogeneous submanifolds. In the context of Teichm\"uller dynamics, the strata of Abelian differentials do not behave like homogeneous spaces. Hence it is unclear whether $\GL_2^{+}(\bbR)$-orbit closures can have nice geometric structures. Nevertheless, the recent breakthrough of 
Eskin-Mirzakhani-Mohammadi (\cite{EskinMirzakhani, EskinMirzakhaniMohammadi}) showed that 
$\GL_2^{+}(\bbR)$-orbit closures are locally cut out by homogeneous linear equations of period coordinates with real coefficients, thus justifying that $\GL_2^{+}(\bbR)$-orbit closures are affine invariant submanifolds.\footnote{Zorich called it the ``magic wand theorem'', which is part of Mirzakhani's Fields Medal work.} Previously this result was only proved in genus two by the fundamental work of McMullen (\cite{McMullenOrbit}). 

Recall that the period coordinates at $(X, \omega) \in \HH(\mu)$ are given by integrating $\omega$ over a basis of the relative homology $H_1(X, \Sigma; \bbZ)$, where $\Sigma$ is the set of zeros of $\omega$. Although period coordinates are not canonical, any two choices of period coordinates differ by a matrix in $\GL_n(\bbZ)$, where $n = \dim \HH(\mu)$, hence the above theorem does not depend on the choice of period coordinates. The proof of the theorem is remarkably long and technical, which involves many ideas from dynamics on homogeneous space, ergodic theory, and measure theory.\footnote{The paper of Eskin-Mirzakhani (\cite{EskinMirzakhani}) is more than $170$ pages. The author heard from Eskin that even a referee report for the paper is more than $40$ pages.} 

Since period coordinates are transcendental, affine invariant submanifolds are apriori only complex-analytic submanifolds. 
The hidden algebraic nature has been discovered by Filip (\cite{Filip}), who proved that affine invariant submanifolds are algebraic subvarieties of $\HH(\mu)$, defined over $\overline{\bbQ}$. In particular, Filip used tools from variations of Hodge structures and showed that affine invariant submanifolds (except those of full rank, 
see Section~\ref{subsec:affine-class}) parameterize curves with non-trivial endormophisms, such as real multiplication on a factor of the Jacobians, which generalizes M\"oller's earlier work (\cite{MoellerHodge, MoellerTorsion}) on torsion and real multiplication for Teichm\"uller curves. As a corollary, the closure of any Teichm\"uller disk (in the standard topology) in $\MM_g$ is a subvariety of $\MM_g$. 

\subsection{Classification of affine invariant submanifolds}
\label{subsec:affine-class}

After the structure theorem of Eskin-Mirzakhani-Mohammadi, the classification of affine invariant submanifolds remains to be a central open question in the study of Teichm\"uller curves. 

The tangent space of an affine invariant submanifold $\calN$ at $(X, \omega)$ can be identified with $H^1(X, \Sigma; \bbC)$ under the period coordinates.\footnote{In Teichm\"uller dynamics $\MM$ is commonly used to denote an affine invariant submanifold, but here we reserve $\MM$ for the moduli space.} Denote by 
$p: H^1(X, \Sigma; \bbC) \to H^1(X; \bbC)$ the projection from the relative cohomology to the absolute cohomology. Avila-Eskin-M\"oller (\cite{AvilaEskinMoeller})
proved that $p(T_{(X, \omega)}\calN)$ is a complex symplectic vector space, hence it has even dimension. 
Define the \emph{rank} of $\calN$ to be $\frac{1}{2} \dim_{\bbC}p(T_	{(X, \omega)}\calN)$, which is at most $g$ by definition. 
If the rank of $\calN$ is bigger than one, we say that $\calN$ is of \emph{higher rank}. If it is equal to $g$, we say that $\calN$
has \emph{full rank}. For example, arithmetic Teichm\"uller curves generated by square-tiled surfaces have rank one. 

By analyzing the boundary of affine invariant submanifolds, Mirzakhani-Wright (\cite{MirzakhaniWright})
proved that if an affine invariant submanifold is of full rank, then it is either a connected component of a stratum or the hyperelliptic locus in a connected component of a stratum. Apisa (\cite{Apisa}) further showed that all affine invariant submanifolds of higher rank in the hyperelliptic strata arise from covering constructions, which gives a coarse classification of affine invariant submanifolds in the hyperelliptic strata modulo finitely many non-arithmetic Teichm\"uller curves and their connected components. Based on previously known examples, Mirzakhani conjectured that higher rank affine invariant submanifolds are either connected components of the strata or arise from covering constructions. Recently counterexamples of Mirzakhani's conjecture are announced by McMullen-Mukamel-Wright and Eskin-McMullen-Mukamel-Wright, whose discoveries rely on special Hurwitz spaces of branched covers, moduli spaces of pointed genus one curves, and a computer search. 

\subsection{Degeneration of Abelian differentials}
\label{subsec:degeneration}

Despite the analytic definition of Teichm\"uller dynamics, a profound algebro-geometric foundation behind the story has already been revealed by many of the preceding  results. In order to apply ideas from algebraic geometry, one upshot is to understanding \emph{degeneration} of Abelian differentials, or equivalently, describing a \emph{compactification} of strata of Abelian differentials, analogous to the Deligne-Mumford compactification $\BM_g$ of the moduli space of curves. 

Recall that the Hodge bundle $\HH$ extends to a rank $g$ bundle $\BHH$ over $\BM_g$, parameterizing stable differentials that are sections of the dualizing line bundle. Hence it would be natural to compactify the stratum $\HH(\mu)$ by taking its closure in $\BHH$. Nevertheless, a disadvantage of this Hodge bundle compactification is that it can lose information of the limit positions of the zeros of $\omega$, especially if $\omega$ vanishes entirely on a component of the underlying reducible curve. Alternatively, up to scaling an Abelian differential is determined by the associated canonical divisor. Hence one can consider the 
\emph{stratum $\calP(\mu)$ of canonical divisors of type $\mu$} as the projectivization of $\HH(\mu)$. Marking the $n$ zeros of the divisors, $\calP(\mu)$ can be regarded as a subvariety of the Deligne-Mumford moduli space $\BM_{g,n}$ of stable genus $g$ curves with $n$ marked points, hence one can study degeneration of canonical divisors of type $\mu$ by taking the closure of $\calP(\mu)$ in $\BM_{g,n}$. 

Eisenbud-Harris (\cite{EisenbudHarrisLimit}) developed a theory of limit linear series that studies degeneration of line bundles and their sections from smooth curves to nodal curves of compact type. In our context the situation is slightly different, because the zero type $\mu$ of the sections is fixed and the underlying curves may fail to be of compact type. Nevertheless, the upshot of \emph{twisting} line bundles by irreducible components of a nodal curve still works. More precisely, define a \emph{twisted canonical divisor of type 
$\mu$} on a nodal curve $X$ to be a collection of (possibly meromorphic) canonical divisors $D_j$ on each irreducible component $X_j$ of $X$ such that the following conditions hold: 
\begin{enumerate}[(1)]
\item The support of $D_j$ is contained in the set of marked points and the nodes lying in $X_j$. Moreover, if $p_i$ is a marked point contained in $X_j$, then $\ord_{p_i}(D_j) = m_i$. 

\item If $q$ is a node of $X$ lying in two irreducible components $X_1$ and $X_2$, then 
$\ord_{q}(D_1) + \ord_{q}(D_2) = -2$. 

\item If $q$ is a node of $X$ lying in two irreducible components $X_1$ and $X_2$ such that 
$\ord_{q}(D_1) = \ord_{q}(D_2) = -1$, then for any node $q' \in X_1\cap X_2$, 
$\ord_{q'}(D_1) = \ord_{q'}(D_2) = -1$. In this case we write $X_1 \sim X_2$. 

\item If $q$ is a node of $X$ lying in two irreducible components $X_1$ and $X_2$ such that 
$\ord_{q}(D_1) > \ord_{q}(D_2)$, then for any node $q' \in X_1\cap X_2$, 
$\ord_{q'}(D_1) > \ord_{q'}(D_2)$. In this case we write $X_1 \succ X_2$. 

\item There does not exist a directed loop $X_1 \succeq X_2 \succeq \cdots \succeq X_k \succeq X_1$ unless all the relations are $\sim$, where $\succeq$ means $\sim$ or $\succ$. 
\end{enumerate}

We briefly explain the motivation behind these conditions. Since the vanishing order along each zero section remains unchanged in a family, it implies condition (1). Since $K_X|_{X_i}$ is locally generated by differentials with a simple pole at a node $q\in X_1\cap X_2$ for $i=1,2$, when twisting by $X_i$, the vanishing order increases by one on one branch of $q$ and decrease by one on the other branch, hence the sum of the vanishing orders does not vary, which implies condition (2). If there is no twist at all, it is the case corresponding to condition (3). Note that twisting by $X_1 + X_2$ does nothing to the nodes lying in their intersection. Hence one can consider twisting, say by a multiple of $X_1$ only. Then the vanishing orders at all nodes between $X_1$ and $X_2$ increase or decrease simultaneously, hence condition (4) follows. By the same token, 
 $X_1 \succeq X_2$ means the twisting coefficient of $X_1$ is bigger than or equal to that of $X_2$, which implies the last condition. 

By an analytic approach, Gendron (\cite{Gendron}) implicitly derived the above conditions and used them to study the Kodaira dimension of strata of twisted canonical divisors. Motivated by the theory of limit linear series, the author (\cite{ChenBoundary})
considered these conditions for curves of compact type and used them to study Weierstrass point behavior for general elements in the strata.  
Farkas and Pandharipande (\cite{FarkasPandharipande}) imposed explicitly these conditions and showed that the corresponding closures in $\BM_{g,n}$ are reducible in general, containing extra boundary components of dimension one less compared to the main component. It is thus natural to ask what extra conditions can distinguish the main component from the other boundary components in the closure. 

\begin{figure}[h]
    \centering
    \psfrag{Y}{$Y_t$}
    \psfrag{Z}{$Z_t$}
    \psfrag{v}{$v_t$}
 \includegraphics[scale=1.0]{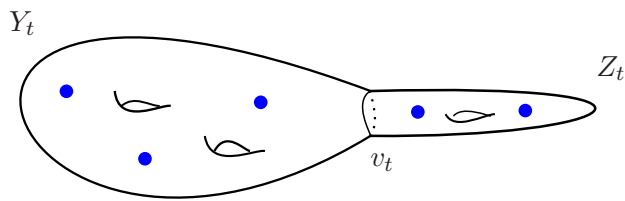}
    \caption{\label{fig:vanishing} A surface with a separating vanishing cycle}
    \end{figure}

Joint with Bainbridge, Gendron, Grushevsky, and M\"oller (\cite{BCGGM1}), we have found the missing condition that arises from a \emph{residue} constraint. Consider the following example. Suppose a family $\calX$ of Abelian differentials $(X_t, \omega_t)$ degenerate to a nodal curve $X$ at $t=0$. Suppose $X$ has a separating node $q$ joining two components $Y$ and $Z$. Without loss of generality, suppose $\lim\limits_{t\to 0} \omega_t |_Y = \eta_Y$ is a holomorphic differential on $Y$, and 
$\lim\limits_{t\to 0} (t^{-\ell}\omega_t) |_Z = \eta_Z$ is a meromorphic differential on $Z$, where $\ell \in \bbZ^{+}$. In other words, the twisting at $q$ is given by $\OO_{\calX}(-\ell Z)$ from the viewpoint of limit linear series. Let $v_t$ be the vanishing cycle on $X_t = Y_t\cup Z_t$ 
that shrinks to the node $q$, where $Y_t\to Y$ and $Z_t\to Z$ as $t\to 0$, see Figure~\ref{fig:vanishing}. 

Since $q$ is a separating node, $v_t = 0 \in H_1(X_t; \bbZ)$ for $t$ nearby $0$. It follows that 
$\int_{v_t} \omega_t = 0$ for $t\neq 0$ and hence $\int_{v_t} t^{-\ell} \omega_t = 0$. Taking the limit as $t\to 0$ and restricting to the component $Z$, we conclude that $\Res_{q} (\eta_Z) = 0$. Conversely, once such a residue condition holds, along with conditions (1)--(4) we are able to prove that the limit differential is smoothable by plumbing techniques in complex-analytic geometry as well as by constructions of flat surfaces. 

\subsection{Cycle classes of strata of Abelian differentials}
\label{subsec:cycle}

As mentioned before, affine invariant submanifolds can provide special subvarieties in the moduli space of curves. It is natural to ask if one can calculate their \emph{cycle classes} in the Chow ring of the moduli space. The first step is to calculate the cycle classes of the strata of Abelian differentials. Tarasca and the author (\cite{ChenTarasca}) calculated the cycle classes of several minimal strata in low genus. Mullane (\cite{Mullane}) obtained a closed formula for the classes of all strata whose projections are effective divisors in $\BM_g$. Both results rely on classical intersection theory on the moduli space of curves. Modulo a conjectural relation to Pixton's formula of the double ramification cycle, 
Janda-Pandharipande-Pixton-Zvonkine (\cite[Appendix]{FarkasPandharipande}) obtained a recursive method to compute the cycle classes of all strata in $\BM_{g,n}$, which suggests that the strata classes are \emph{tautological}. Their approach relies on analyzing the closure of a stratum of twisted canonical divisors in $\BM_{g,n}$ by imposing conditions (1)--(4) in the preceding section. Although the closure may have extra components contained in the boundary, those extra components are products of simpler strata of (possibly meromorphic) differentials, hence one can calculate the class of the main component recursively. 

One can also consider the cycle class calculation in the Hodge bundle compactification. Korotkin-Zograf (\cite{KorotkinZograf}) applied Tau functions to compute the divisor class of the closure of the stratum $\calP(2, 1^{2g-4})$ in the projectivized Hodge bundle over $\BM_g$. The author (\cite{ChenCycle}) gave another proof of this divisor class using intersection theory. Recently Sauvaget-Zvonkine have announced that they are able to compute the cycle classes of all strata closures in the Hodge bundle compactification. 

\section{Meromorphic and higher order differentials}
\label{sec:higher}

In this section we generalize the discussion of Abelian differentials to higher order differentials possibly with poles. 

\subsection{Quadratic differentials}
\label{subsec:quad}

\begin{figure}[h]
    \centering
 \includegraphics[scale=0.8]{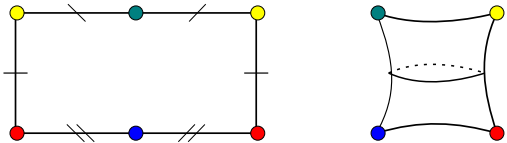}
    \caption{\label{fig:pillow} A quadratic differential with four simple poles on $\bbP^1$}
    \end{figure}

First, if $q$ is a quadratic differential, that is, a section of $K^{\otimes 2}$, then it also induces a flat structure, if one allows \emph{reflection} in addition to translation as transition functions. The flat structure can be defined locally by taking a square root $\omega$ of $q$, which is up to $\pm$, and that is the reason why reflection needs to be part of the transition functions. Moreover if $q$ has at worst simple poles, integrating $\omega$ along a path always provides finite length, hence the corresponding flat surface has finite area. In general, quadratic differentials with at worst simple poles are called \emph{half-translation surfaces}. For example, Figure~\ref{fig:pillow} presents a quadratic differential with four simple poles on $\bbP^1$ as a pillowcase.  

Suppose $(X, q)$ is a half-translation surface. One can take the unique \emph{canonical double cover} $\pi: \hat{X}\to X$ branched at the odd singularities of $q$ (zeros of odd order and simple poles). Then there exists a global Abelian differential $\hat{\omega}$ on $\hat{X}$ such that $\pi^{*}q = \hat{\omega}^2$. From this viewpoint, all questions for Abelian differentials can be similarly asked for quadratic differentials, and indeed many of them can be similarly answered. 

In what follows we mention several results analogous to the case of Abelian differentials. Fix a partition $\nu$ of $4g-4$ such that all entries of $\nu$ are $\geq -1$. Let $\calQ(\nu)$ be the stratum of quadratic differentials of type $\nu$ that are not global squares of Abelian differentials. Lanneau (\cite{LanneauQuad})
classified the connected components of $\calQ(\nu)$ for all $\nu$. In $g\geq 5$, it can have at most two connected components, where extra components are caused by hyperelliptic structures. In particular, the spin parity in the Abelian case does not give rise to additional components in the quadratic case (\cite{LanneauSpin}). 
In genus three and four, there are several \emph{exceptional} disconnected strata. Joint with M\"oller (\cite{ChenMoellerQuadratic})
we found an algebraic parity arising from geometry of canonical curves that distinguishes these exceptional components. We further discovered a number of strata $\calQ(\nu)$ in low genus whose Teichm\"uller curves have non-varying sums of Lyapunov exponents by similar techniques as in the Abelian case. The relation between the area Siegel-Veech constant and the sum of Lyapunov exponents for affine invariant submanifolds in $\calQ(\nu)$ also holds by the seminal work of 
Eskin-Kontsevich-Zorich (\cite{EKZ}). Eskin-Okounkov (\cite{EskinOkounkovQuad}) analyzed the volume growth of $\calQ(\nu)$ by enumerating covers of $\bbP^1$ with certain ramification profile determined by $\nu$. Athreya-Eskin-Zorich (\cite{AthreyaEskinZorich}) proved an explicit formula for the volume of $\calQ(\nu)$ in genus zero. Goujard (\cite{Goujard2}) obtained explicit values for the volume of $\calQ(\nu)$ in all low dimensions. Masur-Zorich (\cite{MasurZorich}) studied the principal boundary of $\calQ(\nu)$, aiming at a recursive way to calculate Siegel-Veech constants of saddle connections. Goujard (\cite{Goujard1}) proved an explicit formula that relates volumes of $\calQ(\nu)$ and Siegel-Veech constants. Grivaux-Hubert (\cite{GrivauxHubert}) constructed explicit affine invariant submanifolds in $\calQ(\nu)$ of arbitrarily large dimension with completely degenerate Lyapunov spectrum. 

\subsection{Differentials with poles}
\label{subsec:pole}

Previously when we talked about an Abelian differential $\omega$, we assumed that $\omega$ is holomorphic, so the corresponding translation surface has finite area. In many cases it would be useful to consider \emph{meromorphic} differentials, in particular when we study the boundary structure of strata of Abelian differentials. 

\begin{figure}[h]
    \centering
    \psfrag{a}{$a$}
     \psfrag{b}{$b$}
      \psfrag{l}{$l_1$}
       \psfrag{m}{$l_2$}
 \includegraphics[scale=0.8]{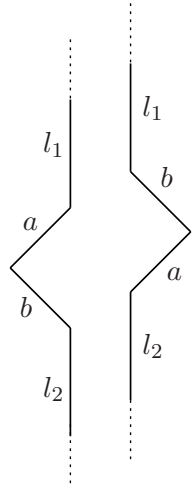}
    \caption{\label{fig:pole-simple} A flat surface with two simple poles}
    \end{figure}

First, we describe the local flat geometry around a pole. Suppose $p$ is a \emph{simple} pole of $\omega$. Then the flat neighborhood of $p$ can be viewed as a \emph{half-infinite cylinder}. The \emph{width} of the cylinder corresponds to the \emph{residue} of $\omega$ at $p$. Figure~\ref{fig:pole-simple} exhibits a meromorphic differential that has two simple poles with opposite residues, where the poles locate at the positive infinity and negative infinity. 

If $p$ is a pole of order $m \geq 2$, Boissy (\cite{Boissy})
showed that one can glue $2m-2$ \emph{broken half-planes} consecutively to form a flat-geometric presentation for a neighborhood of $p$. The boundary of each broken half-plane consists of a half-line to the left and a half-line to the right, which are connected by finitely many broken line segments (saddle connections). The idea behind Boissy's description is the following. If one glues $2m-2$ half-disks consecutively to form a zero of order $m-2$, the local expression of the differential would be $z^{m-2} dz$. Now change coordinates by $z = 1/w$. Then the differential becomes
$\sim w^{-m} dw$, which has a pole of order $m$, and the half-disks turn into half-planes. In particular, the residue of $\omega$ at $p$ is determined by the complex lengths of the boundary line segments and the gluing pattern of the broken half-planes. 

The lower right flat surface of infinite area in Figure~\ref{fig:two-tori} below shows the flat-geometric presentation of a meromorphic differential with a double zero and a double pole on a torus, where the pole locates at infinity. It is constructed by removing the interior of a parallelogram $Z$ from the Euclidean plane and then identifying parallel edges by translation. If we slit the plane along a diagonal of $Z$, we thus recover a pair of broken half-planes as in Boissy's description. In particular, the double pole has no residue.  

Building on earlier work of Kontsevich-Zorich and Lanneau (\cite{KontsevichZorich, LanneauQuad}), Boissy (\cite{Boissy}) classified the connected components of strata of meromorphic differentials with prescribed numbers and multiplicities of zeros and poles. Similarly to the holomorphic case, the strata of meromorphic differentials can have at most three connected components, distinguished by hyperelliptic and spin structures. 

Let us illustrate an interesting viewpoint using flat geometry of meromorphic differentials to study the boundary of strata of Abelian differentials. The flat surface on the left side of Figure~\ref{fig:two-tori} lies in $\HH(2)$, which is constructed by removing a parallelogram $Z$ from the interior of a parallelogram $Y$ and identifying parallel edges. If we shrink $Z$ to a point, we obtain a holomorphic differential $(Y, \eta_Y) \in \HH(0)$, where the marked point encodes the limit position of the inner square. Alternatively, modulo scaling this procedure amounts to expanding $Y$ to be arbitrarily large, hence the limit object represents a meromorphic differential $(Z, \eta_Z) \in \HH(2, -2)$. 
\begin{figure}[h]
    \centering
    \psfrag{Y}{$Y$}
     \psfrag{Z}{$Z$}
      \psfrag{S}{Shrink $Z$}
     \psfrag{E}{Expand $Y$}
 \includegraphics[scale=0.8]{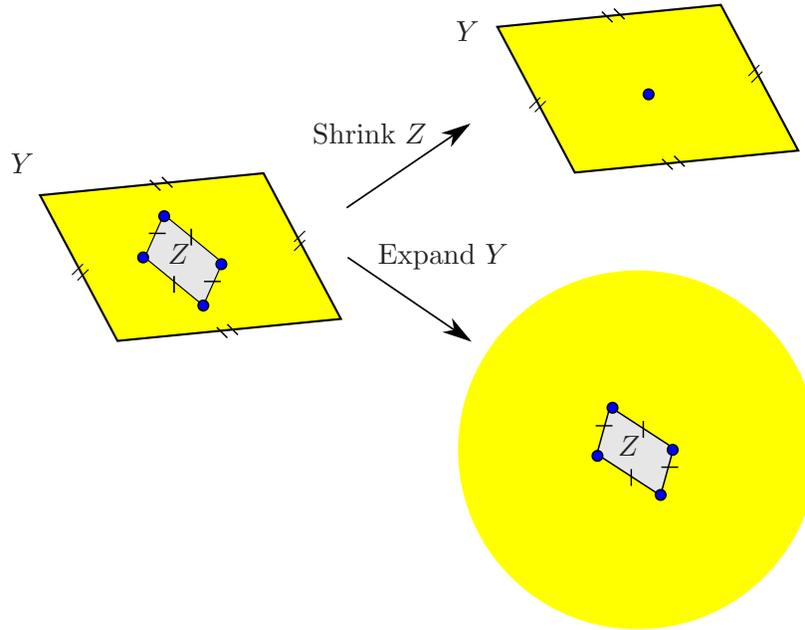}
    \caption{\label{fig:two-tori} Shrinking $Z$ versus expanding $Y$}
    \end{figure}

Note that these two perspectives correspond to exactly the two \emph{aspects} in the theory of limit linear series, which applied to this case says that a family of curves of genus two with a marked double zero of a canonical divisor (Weierstrass point) degenerates to two elliptic curves joined at one node, where the limit marked point is $2$-torsion respect to the node. 

From the viewpoint of $\GL^{+}_2(\bbR)$-action and Teichm\"uller dynamics, the strata of meromorphic differentials somehow display different properties compared to the case of Abelian differentials (see \cite[Appendix]{Boissy}). 

\subsection{Higher order differentials}
\label{subsec:higher}

We conclude the survey by describing some future directions. From the viewpoint of algebraic geometry, an Abelian differential is a section of the canonical line bundle $K$, and a quadratic differential is a section of $K^{\otimes 2}$. Therefore, it is natural to consider \emph{higher order differentials} arising from sections of $K^{\otimes k}$ for a fixed positive integer $k$. Suppose $\mu = (k_1, \ldots, k_n)$ is a partition of $k(2g-2)$, where we allow $k_i$ to be possibly negative. In other words, we want to take meromorphic differentials into account. Let $\HH^k(\mu)$ be the \emph{stratum of $k$-differentials of type $\mu$}, which parameterizes (possibly meromorphic) sections of $K^{\otimes k}$ with zeros and poles of type $\mu$ on genus $g$ Riemann surfaces. All previous questions regarding Abelian and quadratic differentials can be asked similarly for $k$-differentials. In particular, what are the dimension, connected components, compactification, invariants, and cycle class of 
$\HH^k(\mu)$? In a forthcoming work (\cite{BCGGM2}), we will treat these questions systematically. 


\end{document}